\documentclass{amsart}
\usepackage{amsrefs}

\numberwithin{equation}{section}
\newtheorem{thm}{Theorem}
\newtheorem{cor}[thm]{Corollary}
\newtheorem{lemma}[thm]{Lemma}

\newcommand{\abs}[1]{\left\lvert#1\right\rvert}
\newcommand{\dbars}[1]{\left\lVert#1\right\rVert}

\newcommand{\supnorm}[1]{\dbars{#1}_\infty}
\newcommand{\meas}[1]{\dbars{#1}_0}
\newcommand{\Cmplx}[0]{\mathbb{C}}
\newcommand{\Reals}[0]{\mathbb{R}}
\newcommand{\Ints}[0]{\mathbb{Z}}
\newcommand{\AF}[0]{\mathcal{A}}  

\begin{document}

\title{Generalizations of Gon\c{c}alves' inequality}
\author{Peter Borwein}
\address{Department of Mathematics and Statistics\\
Simon Fraser University\\
Burnaby, B.C. V5A~1S6 Canada}
\thanks{Research of P.~Borwein supported in part by NSERC of Canada and MITACS}
\email{pborwein@cecm.sfu.ca}
\author{Michael J. Mossinghoff}
\address{Department of Mathematics\\
Davidson College\\
Davidson, North Carolina 28035 USA}
\email{mjm@member.ams.org}
\author{Jeffrey D. Vaaler}
\address{Department of Mathematics\\
University of Texas\\
Austin, Texas 78712 USA}
\email{vaaler@math.utexas.edu}


\begin{abstract}
Let $F(z)=\sum_{n=0}^N a_n z^n$ be a polynomial with complex coefficients and roots $\alpha_1$, \ldots, $\alpha_N$, let $\dbars{F}_p$ denote its $L_p$ norm over the unit circle, and let $\meas{F}$ denote Mahler's measure of $F$.
Gon\c{c}alves' inequality asserts that
\begin{align*}
\dbars{F}_2
&\geq \abs{a_N} \left( \prod_{n=1}^N \max\{1, \abs{\alpha_n}^2\} + \prod_{n=1}^N \min\{1, \abs{\alpha_n}^2\} \right)^{1/2}\\
&= \meas{F}\left(1+\frac{\abs{a_0 a_N}^2}{\meas{F}^4}\right)^{1/2}.
\end{align*}
We prove that
\[
\dbars{F}_p \geq B_p \abs{a_N} \left( \prod_{n=1}^N \max\{1, \abs{\alpha_n}^p\} + \prod_{n=1}^N \min\{1, \abs{\alpha_n}^p\} \right)^{1/p}
\]
for $1\leq p\leq2$, where $B_p$ is an explicit constant, and that
\[
\dbars{F}_p \geq \meas{F}\left(1+\frac{p^2\abs{a_0 a_N}^2}{4\meas{F}^4}\right)^{1/p}
\]
for $p\geq1$.
We also establish additional lower bounds on the $L_p$ norms of a polynomial in terms of its coefficients.
\end{abstract}

\maketitle

\section{Introduction}\label{sectionIntroduction}

Let $\Delta\subset\Cmplx$ denote the open unit disc, $\overline{\Delta}$ its closure, and let $\AF(\Delta)$ denote the algebra of continuous functions $f :\overline{\Delta}\to\Cmplx$ that are analytic on $\Delta$.
Then $\{\AF(\Delta),\supnorm{\cdot}\}$ is a Banach algebra, where
\[
\supnorm{f}
= \sup\left\{\abs{f(z)} : z\in\overline{\Delta}\right\}
= \sup\left\{\abs{f(e(t))}:t\in\Reals/\Ints\right\},
\]
and $e(t)$ denotes the function $e^{2\pi i t}$.
If $f\in\AF(\Delta)$ and $0<p<\infty$, we also define
\[
\dbars{f}_p = \left(\int_0^1 \abs{f(e(t))}^p dt\right)^{1/p},
\]
and we define
\[
\meas{f} = \exp\left(\int_0^1 \log\abs{f(e(t))} dt\right).
\]
It is known (see \cite{HLP}*{Chapter 6}) that for each $f$ in $\AF(\Delta)$ the function $p\to\dbars{f}_p$ is continuous on $[0,\infty]$, and if $0<p<q<\infty$, then these quantities satisfy the basic inequality
\begin{equation}\label{eqnLpIneqs}
\meas{f} \leq \dbars{f}_p \leq \dbars{f}_q \leq \supnorm{f}.
\end{equation}
Clearly, equality can occur throughout (\ref{eqnLpIneqs}) if $f$ is constant.
On the other hand, if $f$ is not constant in $\AF(\Delta)$, then the function $p\to\dbars{f}_p$ is strictly increasing on $[0,\infty]$.

Now suppose that $F$ is a polynomial in $\Cmplx[z]$ of degree $N\geq1$, and write
\begin{equation}\label{eqnF}
F(z) = \sum_{n=0}^N a_n z^n = a_N \prod_{n=1}^N(z-\alpha_n).
\end{equation}
In this case, the quantity $\meas{F}$ is \textit{Mahler's measure} of $F$, and by Jensen's formula one obtains the well-known identity
\begin{equation}\label{eqnMM}
\meas{F} = \abs{a_N} \prod_{n=1}^N \max\{1, \abs{\alpha_n}\}.
\end{equation}
Thus a special case of (\ref{eqnLpIneqs}) is the inequality (often called Landau's inequality)
\[
\dbars{F}_2 \geq \abs{a_N} \prod_{n=1}^N \max\{1, \abs{\alpha_n}\}.
\]
For polynomials of positive degree, the sharper inequality
\begin{equation}\label{eqnGI1}
\dbars{F}_2 \geq \abs{a_N} \left( \prod_{n=1}^N \max\{1, \abs{\alpha_n}^2\} + \prod_{n=1}^N \min\{1, \abs{\alpha_n}^2\} \right)^{1/2}
\end{equation}
was obtained by Gon\c{c}alves \cite{Goncalves}.
Note that equality occurs in (\ref{eqnGI1}) for constant multiples of $z^N-1$.
Alternatively, the inequality (\ref{eqnGI1}) may be written in the less symmetrical form
\begin{equation}\label{eqnGI2}
\dbars{F}_2 \geq \meas{F} \left( 1 + \frac{\abs{a_0 a_N}^2}{\meas{F}^4}\right)^{1/2}.
\end{equation}

For a positive real number $p$, define the real number $B_p$ by
\begin{equation}\label{eqnBp}
B_p = \left(\frac{1}{2}\int_0^1 \abs{1-e(t)}^p\,dt\right)^{1/p} = \left(\frac{\Gamma(p+1)}{2\Gamma(p/2+1)^2}\right)^{1/p},
\end{equation}
and note that $B_1=2/\pi$ and $B_2=1$.
In this article we establish the following generalizations of Gon\c{c}alves' inequality.

\begin{thm}\label{thmGG}
Let $F(z)\in\Cmplx[z]$ be given by $(\ref{eqnF})$.
If $1\leq p\leq2$, then
\begin{equation}\label{eqnGG1}
\dbars{F}_p \geq B_p \abs{a_N} \left( \prod_{n=1}^N \max\{1, \abs{\alpha_n}^p\} + \prod_{n=1}^N \min\{1, \abs{\alpha_n}^p\} \right)^{1/p}
\end{equation}
and if $p\geq1$, then
\begin{equation}\label{eqnGG2}
\dbars{F}_p \geq \meas{F}\left( 1 + \frac{p^2\abs{a_0 a_N}^2}{4\meas{F}^4}\right)^{1/p}.
\end{equation}
\end{thm}

Equality occurs in (\ref{eqnGG1}) for constant multiples of $z^N-1$.
The inequality (\ref{eqnGG2}) is never sharp for $p\neq2$, but since $B_p<1$ for $1\leq p<2$ it is clearly stronger than (\ref{eqnGG1}) in this range when $\meas{F}^2/\abs{a_0 a_N}$ is large.
For example, one may verify that (\ref{eqnGG2}) produces a better bound in the case $p=1$ whenever
\[
\frac{\meas{F}^2}{\abs{a_0 a_N}} > \frac{\pi^2}{2\left(2-\sqrt{4+2\pi-\pi^2}\right)}= 1.1576382\ldots
\]
Also, for fixed $F$ the right side of (\ref{eqnGG2}) achieves a maximum at $p = 2c\meas{F}^2/\abs{a_0 a_N}$, where $c=1.9802913\ldots$ is the unique positive number satisfying $2c^2 = (1+c^2)\log(1+c^2)$.
In view of (\ref{eqnLpIneqs}), inequality (\ref{eqnGG2}) is therefore only of interest when $1\leq p\leq 2c\meas{F}^2/\abs{a_0 a_N}$.

To prove Theorem~\ref{thmGG}, we first establish some lower bounds on the $L_p$ norms of a polynomial in terms of two of its coefficients $a_L$ and $a_M$, provided $\abs{M-L}$ is sufficiently large.
These inequalities have some independent interest, and we record the results in the following theorem.

\begin{thm}\label{thmLp}
Let $F(z)\in\Cmplx[z]$ be given by $(\ref{eqnF})$, and let $L$ and $M$ be integers satisfying $0\leq L<M\leq N$ and $M-L>\max\{L,N-M\}$.
Then
\begin{equation}\label{eqnBSI}
\supnorm{F} \geq \abs{a_L}+\abs{a_M}.
\end{equation}
Further, if $1\leq p\leq2$ then
\begin{equation}\label{eqnBpI1}
\dbars{F}_p \geq B_p\left(\abs{a_L}^p+\abs{a_M}^p\right)^{1/p},
\end{equation}
and if $p\geq1$ and $a_L$ and $a_M$ are not both~$0$, then
\begin{equation}\label{eqnBpI2}
\dbars{F}_p \geq \max\{\abs{a_L},\abs{a_M}\}\left(1+\left(\frac{p\min\{\abs{a_L},\abs{a_M}\}}{2\max\{\abs{a_L},\abs{a_M}\}}\right)^2\right)^{1/p}.
\end{equation}
\end{thm}

At this point, it is instructive to recall the Hausdorff-Young inequality.
If $p=2$ and $F(z)$ is given by (\ref{eqnF}), then by Parseval's identity we have
\begin{equation}\label{eqnParseval}
\dbars{F}_2 = \left(\abs{a_0}^2+\abs{a_1}^2+\cdots+\abs{a_N}^2\right)^{1/2}.
\end{equation}
If $p=1$, then the inequality
\begin{equation}\label{eqnEasyL1}
\dbars{F}_1 \geq \max\{\abs{a_0}, \abs{a_1}, \ldots, \abs{a_N}\}
\end{equation}
follows immediately from the identity
\[
a_n = \int_0^1 F(e(t)) e(-nt)\,dt.
\]
Now suppose that $1<p<2$ and let $q$ be the conjugate exponent for $p$, so $p^{-1}+q^{-1}=1$.
Then the Hausdorff-Young inequality \cite{Katz04}*{p.\ 123} asserts that
\begin{equation}\label{eqnHYI}
\dbars{F}_p \geq \left(\abs{a_0}^q+\abs{a_1}^q+\cdots+\abs{a_N}^q\right)^{1/q},
\end{equation}
and so interpolates between (\ref{eqnParseval}) and (\ref{eqnEasyL1}).
If $p=2$, then (\ref{eqnBpI1}) and (\ref{eqnBpI2}) are equivalent and clearly follow from the identity (\ref{eqnParseval}).
But for $1<p<2$, the inequalities (\ref{eqnBpI1}) and (\ref{eqnBpI2}) are not immediate consequences of (\ref{eqnHYI}).
In fact, it is easy to see that the lower bounds in (\ref{eqnBpI1}), (\ref{eqnBpI2}), and (\ref{eqnHYI}) are not comparable.
If $p=1$, the same remarks apply to (\ref{eqnBpI1}), (\ref{eqnBpI2}), and (\ref{eqnEasyL1}).

In section~\ref{sectionBinomials} we develop some preliminary results concerning lower bounds on $L_p$ norms of binomials, and we use these facts to establish Theorems~\ref{thmGG} and~\ref{thmLp} in section~\ref{sectionProofs}.

\section{Norms of binomials}\label{sectionBinomials}

For $0<r<1$ and real $t$, recall that the Poisson kernel is defined by
\[
P(r,t) = \sum_{n=-\infty}^\infty r^{\abs{n}} e(nt) = \Re\left(\frac{1+re(t)}{1-re(t)}\right) = \frac{1-r^2}{\abs{1-re(t)}^2}.
\]
This is a positive summability kernel that satisfies
\[
\int_0^1 P(r,t)\,dt = 1
\]
and
\[
\lim_{r\to1-} \int_\epsilon^{1-\epsilon} P(r,t)\, dt = 0
\]
for $0<\epsilon<1/2$.

\begin{lemma}\label{lemPoisson}
If $p>0$ then
\[
\lim_{r\to1-} \int_0^1 \abs{1-re(t)}^p P(r,t)\,dt = 0.
\]
\end{lemma}

\begin{proof}
Let $0<\epsilon<1/2$ so that
\begin{align*}
\int_0^1 \abs{1-re(t)}^p P(r,t)\,dt
&= \int_{-\epsilon}^\epsilon \abs{1-re(t)}^p P(r,t)\,dt + \int_\epsilon^{1-\epsilon} \abs{1-re(t)}^p P(r,t)\,dt\\
&\leq \abs{1-re(\epsilon)}^p \int_{-\epsilon}^\epsilon P(r,t)\,dt + 2^p \int_\epsilon^{1-\epsilon} P(r,t)\,dt\\
&\leq \abs{1-re(\epsilon)}^p + 2^p \int_\epsilon^{1-\epsilon} P(r,t)\,dt.
\end{align*}
We conclude that
\[
\limsup_{r\to1-} \int_0^1 \abs{1-re(t)}^p P(r,t)\,dt \leq \abs{1-e(\epsilon)}^p \leq (2\pi\epsilon)^p,
\]
and the statement follows.
\end{proof}

For positive numbers $p$ and $r$, we define
\begin{equation}\label{eqnIdefn}
\mathcal{I}_p(r) = \int_0^1 \abs{1-r^{1/p}e(t)}^p\,dt.
\end{equation}
It follows easily that $r\to\mathcal{I}_p(r)$ is a continuous, positive, real-valued function that satisfies the functional equation
\begin{equation}\label{eqnIFcnEqn}
\mathcal{I}_p(r) = r \mathcal{I}_p(1/r)
\end{equation}
for all positive $r$.
The following lemma records some further information about this function.

\begin{lemma}\label{lemIinfo}
For any positive number $p$, the function $r\to\mathcal{I}_p(r)$ has a continuous derivative at each point of $(0,\infty)$ and satisfies the identity $\mathcal{I}'_p(1) = \mathcal{I}_p(1)/2$.
Moreover, this function has infinitely many continuous derivatives on the open subintervals $(0,1)$ and $(1,\infty)$.
\end{lemma}

\begin{proof}
Suppose first that $0<r<1$.
Then $\left(1-r^{1/p}e(t)\right)^{p/2}$ has the absolutely convergent Fourier expansion
\[
\left(1-r^{1/p}e(t)\right)^{p/2} = \sum_{m\geq0} \binom{p/2}{m} (-1)^mr^{m/p} e(mt).
\]
By Parseval's identity, we have
\begin{equation}\label{eqnIexpn}
\mathcal{I}_p(r) = \sum_{m\geq0} \binom{p/2}{m}^2 r^{2m/p}.
\end{equation}
This shows that $r\to\mathcal{I}_p(r)$ is represented on $(0,1)$ by a convergent power series in $r^{1/p}$ and therefore has infinitely many continuous derivatives on this interval.
Next, we observe that
\[
\frac{\partial}{\partial r} \abs{1-r^{1/p} e(t)}^p = \frac{\abs{1-r^{1/p} e(t)}^p}{2r}\left(1-P(r^{1/p}, t)\right).
\]
It follows that if $0<\epsilon\leq1/4$ and $\epsilon\leq r\leq1-\epsilon$, then there exists a positive constant $C(\epsilon,p)$  such that
\[
\abs{\frac{\partial}{\partial r} \abs{1-r^{1/p} e(t)}^p} \leq C(\epsilon,p).
\]
From the mean value theorem and the dominated convergence theorem, we find that
\[
\mathcal{I}'_p(r) = \frac{1}{2r}\int_0^1 \abs{1-r^{1/p} e(t)}^p \left(1-P(r^{1/p},t)\right)\,dt,
\]
and therefore
\[
\mathcal{I}_p(r) - 2r\mathcal{I}'_p(r) = \int_0^1 \abs{1-r^{1/p} e(t)}^p P(r^{1/p},t)\,dt.
\]
Using the continuity of $r\to\mathcal{I}_p(r)$ and Lemma~\ref{lemPoisson}, we conclude that
\begin{equation}\label{eqnLeftDeriv}
\lim_{r\to1-} \mathcal{I}'_p(r) = \mathcal{I}_p(1)/2.
\end{equation}
Again using the mean value theorem, it follows that $r\to\mathcal{I}_p(r)$ has a left-hand derivative at~1 with the value $\mathcal{I}_p(1)/2$.

From (\ref{eqnIFcnEqn}) and (\ref{eqnIexpn}) we find that
\begin{equation}\label{eqnIExpnBigr}
\mathcal{I}_p(r) = r\left(\sum_{m\geq0} \binom{p/2}{m}^2 r^{-2m/p}\right)
\end{equation}
for $r>1$.
Thus $r\to\mathcal{I}_p(r)$ is represented by $r$ times a convergent power series in $r^{-1/p}$, and so has infinitely many continuous derivatives on the interval $(1,\infty)$.
Next, we differentiate both sides of (\ref{eqnIFcnEqn}) to obtain the identity
\[
\mathcal{I}'_p(r) = \mathcal{I}_p(1/r) - \frac{\mathcal{I}'_p(1/r)}{r}
\]
for $r>1$, and using the continuity of $r\to\mathcal{I}_p(r)$ and (\ref{eqnLeftDeriv}), we conclude that
\[
\lim_{r\to1+} \mathcal{I}'_p(r) = \mathcal{I}_p(1) - \lim_{s\to1-} \mathcal{I}'_p(s) = \mathcal{I}_p(1)/2.
\]
It follows that $r\to\mathcal{I}_p(r)$ has a right-hand derivative at~1 with value $\mathcal{I}_p(1)/2$.

We conclude then that $r\to\mathcal{I}_p(r)$ is continuously differentiable on $(0,\infty)$ and $\mathcal{I}'_p(1)=\mathcal{I}_p(1)/2$.
\end{proof}

From the proof of the lemma we obtain the following lower bound on the $L_p$ norm of a binomial.

\begin{cor}\label{corBinomialBound1}
Let $0\leq L<M$ be integers and let $\alpha$ and $\beta$ be complex numbers, not both zero.
If $p>0$ then
\[
\dbars{\alpha z^L+\beta z^M}_p \geq\max\{\abs{\alpha},\abs{\beta}\}\left(1+\left(\frac{p\min\{\abs{\alpha},\abs{\beta}\}}{2\max\{\abs{\alpha},\abs{\beta}\}}\right)^2\right)^{1/p},
\]
with equality precisely when $\alpha\beta=0$ or $p=2$.
\end{cor}

\begin{proof}
The result is trivial if either $\alpha$ or $\beta$ is zero, so we assume that this is not the case.
We may then assume by homogeneity that $\alpha=1$, and
it is clear from the definition of $\dbars{f}_p$ that we may assume that $L=0$, and that $\beta$ is real and negative.
If $\abs{\beta}<1$, then taking $r=\abs{\beta}^p$ in (\ref{eqnIexpn}) and keeping just the first two terms of the sum, we obtain
\[
\dbars{1+\beta z^M}_p^p  = \dbars{1+\beta z}_p^p \geq 1 + \frac{p^2\abs{\beta}^2}{4}.
\]
If $\abs{\beta}>1$, then
\[
\dbars{1+\beta z^M}_p^p 
= \abs{\beta}^p \dbars{\beta^{-1}+z}_p^p
= \abs{\beta}^p \dbars{1+z/\beta}_p^p,
\]
so taking $r=\abs{\beta}^{-p}$, we obtain in the same way
\[
\dbars{1+\beta z^M}_p^p \geq \abs{\beta}^p \left(1 + \frac{p^2}{4\abs{\beta}^2}\right).
\]
The case $\beta=-1$ follows by continuity.
For the case of equality, notice that the sum (\ref{eqnIexpn}) has precisely two nonzero terms only when $p=2$.
\end{proof}

The next lower bound is obtained by establishing the convexity of the function $r\to\mathcal{I}_p(r)$ for each fixed $p$ in $(0,2]$.

\begin{lemma}\label{lemIBound}
If $0<p\leq2$, then the function $r\to\mathcal{I}_p(r)$ satisfies the inequality
\begin{equation}\label{eqnIBound}
\mathcal{I}_p(r) \geq \frac{\mathcal{I}_p(1)(1+r)}{2}
\end{equation}
for $r>0$.
\end{lemma}

\begin{proof}
If $p=2$ then $\mathcal{I}_2(r)=1+r$ and the result is trivial.
Suppose then that $0<p<2$.
If $r<1$, then we may differentiate the power series (\ref{eqnIexpn}) termwise to obtain
\[
\mathcal{I}'_p(r) = \sum_{m\geq0} \binom{p/2}{m}^2 \frac{2m}{p} r^{(2m/p)-1}.
\]
As $0<p<2$, it follows that $r\to\mathcal{I}'_p(r)$ is strictly increasing on $(0,1)$, so $r\to\mathcal{I}_p(r)$ is strictly convex on this interval.
Thus, if $r$ and $s$ are in $(0,1)$, then
\begin{equation}\label{eqnConvexity}
\mathcal{I}_p(r) \geq \mathcal{I}_p(s) + (r-s)\mathcal{I}'_p(s).
\end{equation}
Letting $s\to1-$ and using Lemma~\ref{lemIinfo}, we obtain
\begin{equation}\label{eqnConvexity2}
\mathcal{I}'_p(r) \geq \mathcal{I}_p(1) + \mathcal{I}'_p(1)(r-1) = \frac{\mathcal{I}_p(1)(1+r)}{2}.
\end{equation}
for $0<r<1$.

In a similar manner, if $r>1$ we differentiate (\ref{eqnIExpnBigr}) termwise to obtain
\[
\mathcal{I}'_p(r) = \sum_{m\geq0} \binom{p/2}{m}^2\left(1-\frac{2m}{p}\right)r^{-2m/p},
\]
and again $r\to\mathcal{I}'_p(r)$ is strictly increasing on $(1,\infty)$, so $r\to\mathcal{I}_p(r)$ is strictly convex on this interval.
Thus (\ref{eqnConvexity}) holds as well for $r>1$ and $s>1$, and letting $s\to1+$ we obtain (\ref{eqnConvexity2}) for $r>1$.

We have therefore verified (\ref{eqnIBound}) at each point $r$ in $(0,1)\cup(1,\infty)$, and it is trivial at $r=1$.
\end{proof}

Using this lemma, we obtain a second lower bound on the $L_p$ norm of a binomial.

\begin{cor}\label{corBinomialBound2}
Let $0\leq L< M$ be integers and let $\alpha$ and $\beta$ be complex numbers.
If $0<p\leq2$ then
\[
\dbars{\alpha z^L + \beta z^M}_p \geq B_p \left(\abs{\alpha}^p+\abs{\beta}^p\right)^{1/p}.
\]
\end{cor}

\begin{proof}
The result is trivial if either $\alpha$ or $\beta$ is zero, so we assume that this is not the case.
By homogeneity, we may assume then that $\abs{\alpha}=1$, and we may assume that $L=0$ and that $\beta$ is real and negative by the definition of $\dbars{f}_p$.
Using Lemma~\ref{lemIBound}, we obtain
\begin{align*}
\dbars{1+\beta z^M}_p^p
&= \dbars{1+\beta z}_p^p\\
&= \mathcal{I}_p(\abs{\beta}^p)\\
&\geq \frac{\mathcal{I}_p(1)\left(1+\abs{\beta}^p\right)}{2}\\
&= B_p^p\left(1+\abs{\beta}^p\right).
\end{align*}
\end{proof}

\section{Proofs of the theorems}\label{sectionProofs}

The proof of Theorem~\ref{thmLp} employs an averaging argument and makes use of the triangle inequality for $L_p$ norms.
We therefore require the restriction $p\geq1$ in the statement of the theorem.

\begin{proof}[Proof of Theorem~\ref{thmLp}]
Suppose that $F(z)=\sum_{n=0}^N a_n z^n$ is a polynomial with complex coefficients, $L$ and $M$ are as in the statement of the theorem, and $1\leq p\leq2$.
Set $K=M-L$, and let $\zeta_K$ denote a primitive $K$th root of unity in $\Cmplx$.
Then
\begin{align*}
\frac{1}{K} \sum_{k=1}^K \zeta_K^{-kL} F\left(\zeta_K^k z\right)
&= \frac{1}{K} \sum_{n=0}^N \left(\sum_{k=1}^K \zeta_K^{k(n-L)}\right) a_n z^n\\
&= \sum_{\substack{0\leq n\leq N\\ n\equiv L \textrm{\ (mod $K$)}}} a_n z^n\\
&= a_L z^L + a_M z^M.
\end{align*}
Using the triangle inequality and the fact that the polynomials $\zeta_K^{-kL} F(\zeta_K^k z)$ all have the same $L_p$ norm, we find that
\begin{equation}\label{eqnPolyNormIneq}
\dbars{F}_p \geq \dbars{\frac{1}{K} \sum_{k=1}^K \zeta_K^{-kL} F\left(\zeta_K^k z\right)}_p = \dbars{a_L z^L + a_M z^M}_p
\end{equation}
for $1\leq p\leq\infty$.
The inequality (\ref{eqnBSI}) then follows by selecting a complex number $z$ of unit modulus so that $a_L z^L$ and $a_M z^M$ have the same argument.
Then inequalities (\ref{eqnBpI1}) and (\ref{eqnBpI2}) are established by combining (\ref{eqnPolyNormIneq}) with Corollary~\ref{corBinomialBound1} and Corollary~\ref{corBinomialBound2}, respectively.
\end{proof}

The proof of Theorem~\ref{thmGG} proceeds by applying Theorem~\ref{thmLp} to a polynomial having the same values over the unit circle as the given polynomial $F$.
Ostrowski \cite{Ostrowski} and  Mignotte \cite{Mignotte} (see also \cite{MS}*{p.\ 80}) employ a similar construction in their proofs of Gon\c{c}alves' inequality (\ref{eqnGI2}) in the case $p=2$.

\begin{proof}[Proof of Theorem~\ref{thmGG}]
Suppose that $F(z)=\sum_{n=0}^N a_n z^n = a_N\prod_{n=1}^N (z-\alpha_n)$ is a polynomial with complex coefficients.
If $F(z)$ has a root at $z=0$, then (\ref{eqnGG1}) and (\ref{eqnGG2}) follow immediately from (\ref{eqnLpIneqs}), so we assume that $a_0\neq0$.
Let $\mathcal{E}$ denote the collection of all subsets of $\{1,2,\ldots,N\}$, and for each $E$ in $\mathcal{E}$, let $E'$ denote the complement of $E$ in $\{1,2,\ldots,N\}$.
For each set $E$ in $\mathcal{E}$, we define the finite Blaschke product $B_E(z)$ by
\[
B_E(z) = \prod_{n\in E} \frac{1-\overline{\alpha_n}z}{z-\alpha_n}
\]
and the polynomial $G_E(z)$ by
\[
G_E(z) = B_E(z) F(z) = \sum_{n=0}^N b_n(E) z^n.
\]
Clearly,
\[
b_0(E) = a_N \prod_{m\in E'} (-\alpha_m)
\]
and
\[
b_N(E) = a_N \prod_{n\in E} (-\overline{\alpha_n}).
\]
If $\abs{z}=1$ then the Blaschke product satisfies $\abs{B_E(z)}=1$, so
$\dbars{G_E(z)}_p = \dbars{F}_p$ for $0\leq p\leq\infty$ and every $E$ in $\mathcal{E}$.
Now select $L=0$ and $M=N$ for the polynomial $G_E(z)$ in Theorem~\ref{thmLp}.
Then from (\ref{eqnBpI1}) we obtain
\begin{equation}\label{eqnGenGG1}
\dbars{F}_p \geq B_p \abs{a_N}\left(\prod_{m\in E'} \abs{\alpha_m}^p + \prod_{n\in E} \abs{\alpha_n}^p \right)^{1/p}
\end{equation}
for $1\leq p\leq2$.
Also, assuming without loss of generality that $\abs{b_0(E)}\geq\abs{b_N(E)}$, we find from (\ref{eqnBpI2}) that 
\begin{equation}\label{eqnGenGG2}
\begin{split}
\dbars{F}_p
&\geq \abs{b_0(E)}\left(1+\frac{p^2\abs{b_N(E)}^2}{4\abs{b_0(E)}^2}\right)^{1/p}\\
&= \abs{b_0(E)}\left(1+\frac{p^2\abs{a_0 a_N}^2}{4\abs{b_0(E)}^4}\right)^{1/p}
\end{split}
\end{equation}
for $p\geq1$.
Inequalities (\ref{eqnGG1}) and (\ref{eqnGG2}) then follow from (\ref{eqnGenGG1}) and (\ref{eqnGenGG2}) by choosing $E = \{n : \abs{\alpha_n} \leq 1\}$.
\end{proof}

We remark that the choice of $E$ in the preceding proof produces the best possible inequality in (\ref{eqnGenGG1}).
To establish this, suppose that $E\in\mathcal{E}$ has $\abs{b_0(E)}=\meas{F}/r$ for some real number $r$, so $\abs{b_N(E)}=r\abs{a_0 a_N}/\meas{F}$.
Then certainly $1\leq r\leq\meas{F}^2/\abs{a_0 a_N}$, and it is easy to check that
\[
\frac{\meas{F}}{r} + \frac{r\abs{a_0 a_N}}{\meas{F}} \leq \meas{F} + \frac{\abs{a_0 a_N}}{\meas{F}}
\]
in this range, with equality occurring only at the endpoints.

It is possible, however, that a different choice for $E$ in (\ref{eqnGenGG2}) could produce a bound better than (\ref{eqnGG2}) for a particular polynomial.
Specifically, if $E\in\mathcal{E}$ has $\abs{b_0(E)}=\meas{F}/r$, again with $1\leq r\leq\meas{F}^2/\abs{a_0 a_N}$, then we obtain an improved bound whenever
\[
\left(4\meas{F}^4 + p^2\abs{a_0 a_N}^2\right)r^p <  4\meas{F}^4 + p^2\abs{a_0 a_N}^2r^4,
\]
and this may occur when $p$ is small.
For example, the polynomial $F(z)=18z^2-101z+90$ has roots $\alpha_1=9/2$ and $\alpha_2=10/9$; choosing $E=\{\}$ with $p=1$ yields $\dbars{F}_1\geq 90.9$, but selecting $E=\{1\}$ (so $r=9/2$) produces a lower bound slightly larger than 102.


\end{document}